 
 
 
 

  \ifx\osumriploaded\MYundefined  
  \else
    \immediate\write16{osumrip.sty ALREADY loaded.}
  \endinput\fi

  \def\osumriploaded{\relax}

  \magnification=1200 
  \hoffset=.35truein  
  \hsize=355pt 
  \vsize = 500 pt
  \baselineskip=11pt 
  \lineskip=1.1pt
  \lineskiplimit=.8pt
  \topskip=12pt 
  \bigskipamount=10pt plus 4pt minus 1pt
  \hfuzz=\hsize
  \overfullrule = 0 pt 
  \vbadness=10000
  \hbadness=10000
  \nopagenumbers
  \widowpenalty=5000
  \pageno=1
  \footline{\ifnum \pageno>0 \hss\tenrm\folio\hss\fi}
  \mathsurround=1pt
  \def\prose{\kern\mathsurround}

  \def\StdPretolerance{100}
  \tolerance=\StdPretolerance

  \def\StdParskip{0pt}   
  \parskip=\StdParskip
  \parindent=0.5cm
 
  \def\tenpoint{}
 
 
  \font\Bigbf=cmbx10 at 14.4pt 
  \font\twelvebf=cmbx12
  \font\tenbf=cmbx10
  
  \font\tensmc=cmcsc10

 \def\smc{\tensmc}
 
 \def \Smallfonts {\relax } 
  \ifx\amspptloaded@AmS\relax
     \def\Smallfonts{\eightpoint}
  \fi
  \def\next{AMSPPT}
  \ifx\styname\next 
     \def\Smallfonts{\eightpoint}
  \fi

  \def\Hfont{\Bigbf}
  \def\Authorfont{\twelvebf}
  \def\HHfont{\twelvebf} 
  \def\HHHfont{\tenbf}
  \def\Figurefont{\bf}
  \def\Bibfont{\tenbf}
  
 \def \thfont {\smc }
 \def \pffont {\smc }
 \def \rkfont {\smc }
 \def \dffont {\smc }
 \def \egfont {\smc }
 

 
  \chardef\CatAt\the\catcode`\@
  \catcode`\@=11

 \let\wlog@ld\wlog 
 \def\wlog#1{\relax}

   \def\hexnumber@#1{\ifcase#1 
     0\or1\or2\or3\or4\or5\or6\or7\or8\or9\or 
     A\or B\or C\or D\or E\or F\fi}

 
 \newif\ifIN@
 \def\m@rker{\m@@rker}
 \def\IN@{\expandafter\INN@\expandafter}
 \long\def\INN@0#1@#2@{\long\def\NI@##1#1##2##3\ENDNI@
    {\ifx\m@rker##2\IN@false\else\IN@true\fi}%
     \expandafter\NI@#2@@#1\m@rker\ENDNI@}

  \newtoks\Initialtoks@  \newtoks\Terminaltoks@
  \def\SPLIT@{\expandafter\SPLITT@\expandafter}
  \def\SPLITT@0#1@#2@{\def\TTILPS@##1#1##2@{%
     \Initialtoks@{##1}\Terminaltoks@{##2}}\expandafter\TTILPS@#2@}


  \newtoks\Trimtoks@

 \def\ForeTrim@{\expandafter\ForeTrim@@\expandafter}
 \def\ForePrim@0 #1@{\Trimtoks@{#1}}
 \def\ForeTrim@@0#1@{\IN@0\m@rker. @\m@rker.#1@%
     \ifIN@\ForePrim@0#1@%
     \else\Trimtoks@\expandafter{#1}\fi}

  \def\Trim@0#1@{%
      \ForeTrim@0#1@%
      \IN@0 @\the\Trimtoks@ @%
        \ifIN@ 
             \SPLIT@0 @\the\Trimtoks@ @\Trimtoks@\Initialtoks@
             \IN@0\the\Terminaltoks@ @ @%
                 \ifIN@
                 \else \Trimtoks@ {FigNameWithSpace}%
                 \fi
        \fi
      }

 
 
 \def\Hrule{\hrule width0pt height0pt}

 \newskip\LastSkip
 \def\SaveLastSkip{\LastSkip\lastskip}

 \def\NoindentAfter{\everypar={\setbox0=\lastbox\everypar={}}}

 \long\def\H#1\par#2\par{\notenumber=0%
    \hbox to 0pt{}\vskip15pt plus 1pt 
    {\baselineskip=15pt plus 1pt\parindent=0pt\parskip=0pt\frenchspacing
    \leftskip=0pt plus .2\hsize minus .1\hsize 
    \rightskip=0pt plus .2\hsize minus .1\hsize 
    \def\\{\unskip\break} 
    \pretolerance=10000 \Hfont #1\unskip\break
    \vskip20pt plus 1pt\Authorfont #2\unskip\break\par}
    \vskip20pt plus 1pt%
    \NoindentAfter\par\rm
    }




 \newdimen\PageRemainder
  \def\SetPageRemainder{\PageRemainder=\pagegoal\advance\PageRemainder by
   -1\pagetotal}
 
  \def\Rpt@{}\def\Rpt@@{}
  
  \long\def\HH#1\par{\par
  \SaveLastSkip\removelastskip\goodbreak
  \ifdim\LastSkip<24pt
     \LastSkip 24pt plus 1pt\fi
  \SetPageRemainder\advance\PageRemainder-\LastSkip
  \ifdim\PageRemainder<150pt
       \edef\Rpt@{remain = \the\PageRemainder\noexpand\\
                pagetotal=\the\pagetotal\noexpand\\
                           pagegoal=\the\pagegoal}%
          \else\def\Rpt@{}\fi
   \ifdim\PageRemainder<72pt 
             \edef\Rpt@@{\noexpand\\
                      Had HH PageRemainder$<$72pt\noexpand\\
                      Hence forced break!}%
     \vskip 0pt plus .1\PageRemainder\eject 
    \else\edef\Rpt@@{}
    \fi
    \vskip\LastSkip\Hrule 
    \pretolerance=10000\rightskip=0pt plus 3em
    \hangafter1 \hangindent=2.2em%
    \noindent
    \HHfont \unskip \Ednote{\Rpt@\Rpt@@}%
            \def\Rpt@{}\def\Rpt@@{}%
            \ignorespaces
            #1\par\rightskip=0pt\pretolerance=\StdPretolerance%
      \NoindentAfter\tenpoint\rm\vskip 8pt plus 1pt minus 1pt}
 
  \long\def\HHH#1\par{\par%
  \SaveLastSkip\removelastskip\goodbreak
  \ifdim\LastSkip<10pt
     \LastSkip 10pt plus 1pt\fi
  \SetPageRemainder\advance\PageRemainder-\LastSkip
  \ifdim\PageRemainder<150pt
       \edef\Rpt@{remain = \the\PageRemainder\noexpand\\
                pagetotal=\the\pagetotal\noexpand\\
                           pagegoal=\the\pagegoal}%
          \else\def\Rpt@{}\fi
   \ifdim\PageRemainder<60pt  
             \edef\Rpt@@{\noexpand\\
                      Had HHH PageRemainder$<$60pt\noexpand\\
                      Hence forced break!}%
       \vskip 0pt plus .1\PageRemainder\eject 
   \else\edef\Rpt@@{}
   \fi
   \vskip\LastSkip\Hrule\par\noindent
   \HHHfont \unskip\Ednote{\Rpt@\Rpt@@}\ignorespaces
   #1\unskip.\quad\rm\ignorespaces
   \ignorepars}
       
  \long\def\ignorepars#1\par{\def\Test{#1}%
     \ifx\Test\Empty\def\This{\ignorepars}%
        \else\def\This{\Test\par}\fi
           \This}

 
 \def\ProcBreak{\par\ifdim\lastskip<8pt
    \removelastskip
    \penalty-200\vskip8pt plus1pt\fi}

 \def \th #1\par{\ProcBreak \noindent
   {\thfont\ignorespaces #1\unskip.}\quad\it}
 \def \endth {\ProcBreak\rm }
 
 
 \def\pf #1\par{\ProcBreak %
    \noindent\pffont#1\unskip:\kern.5pt
       \raise-.4pt\hbox{---}\quad\rm}
 

  \def\qedbox{\hbox{\vbox{
    \hrule width0.2cm height0.2pt 
    \hbox to 0.2cm{\vrule height 0.2cm width 0.2pt 
             \hfil\vrule height0.2cm width 0.2pt}
    \hrule width0.2cm height 0.2pt}}}

  \def\qed{\ifmmode\qedbox
    \else\hglue4mm\unskip\hfill\qedbox\ProcBreak\fi} 

  \def \rk #1\par{\ProcBreak
     \noindent{\rkfont\ignorespaces #1\unskip.}\quad\rm}

  \def \df #1\par{\ProcBreak 
     \noindent{\dffont\ignorespaces #1\unskip.}\quad\rm}

  \def \eg #1\par{\ProcBreak 
     \noindent{\egfont\ignorespaces #1\unskip.}\quad\rm}



  \newdimen\Overhang

   \def\MaxTag@#1#2#3#4#5{\setbox0=\hbox{#4\ignorespaces#2\unskip}%
     \dimen0=\wd0\advance\dimen0 by#3
     \ifdim\dimen0<#5\relax\dimen0=#5\fi
     \expandafter\edef\csname #1Hang\endcsname{\the\dimen0}}

 \def\MaxItemTag#1{\MaxTag@{Item}{#1}{.4em}{\ItemStyle}{\parindent}}%
 \def\MaxItemItemTag#1{%
        \MaxTag@{ItemItem}{#1}{.4em}{\ItemItemStyle}{\parindent}}
 \def\MaxNrTag#1{\MaxTag@{Nr}{#1}{.5em}{\NrStyle}{\parindent}}
 \def\MaxReferenceTag#1{%
        \MaxTag@{Reference}{[#1]}{.6em}{\Smallfonts}{\parindent}}
 \def\MaxFootTag#1{\MaxTag@{Foot}{#1}{.4em}{\Smallfonts}{\z@}}

  \def\SetOverhang@{\Overhang=.8\dimen0%
     \advance\Overhang by \wd0\relax
     \ifdim\Overhang>\hangindent\relax
       \advance\Overhang by .25\dimen0%
       \Ednote{Tag is pushing text.}
     \else\Overhang=\hangindent
     \fi}

   \def\Item#1{\par\noindent
      \hangafter1\hangindent=\ItemHang
      \setbox0=\hbox{\ItemStyle\ignorespaces#1\unskip}%
      \dimen0=.4em\SetOverhang@
      \rlap{\box0}\kern\Overhang\ignorespaces}

   \def\ItemStyle{\rm}
   \MaxItemTag{(iii)}
         
   \def\ItemItem#1{\par\noindent
      \hangafter1\hangindent=\ItemItemHang
      \setbox0=\hbox{\ItemItemStyle\ignorespaces#1\unskip}%
      \dimen0=.4em\SetOverhang@
      \advance\hangindent by \ItemHang
      \kern\ItemHang\rlap{\box0}%
      \kern\Overhang\ignorespaces}

    \def\ItemItemStyle{\rm}
    \MaxItemItemTag{(iii)}

  \def\Nr#1{\par\noindent\hangindent=\NrHang 
    \setbox0=\hbox{\NrStyle\ignorespaces#1\unskip}%
    \dimen0=.5em\SetOverhang@
    \rlap{\box0}\kern\Overhang
    \hangindent=\z@\ignorespaces}

   \def\NrStyle{\rm}
   \MaxNrTag{(2)}

   \newskip\Rosterskip\Rosterskip 1pt plus1pt 
   \def\Roster{\par\ifdim\lastskip<\Rosterskip\removelastskip\vskip\Rosterskip\fi
    \bgroup}
   \def\endRoster{\par\global\edef\LastSkip@{\the\lastskip}\removelastskip
       \egroup\penalty-50\LastSkip\LastSkip@\relax
       \ifdim\LastSkip<\Rosterskip\LastSkip\Rosterskip\fi
       \vskip\LastSkip}




\def\cite#1{
    \def\nextiii@##1,##2\end@{{\frenchspacing\bf
      \lBr\ignorespaces##1\unskip{\rm,~\ignorespaces##2}\rBr}}%
    \IN@0,@#1@%
    \ifIN@\def\next{\nextiii@#1\end@}\else
    \def\next{{\bf\lBr#1\rBr}}\fi\next}



  
   \def \Bib#1\par{%
       \par\removelastskip\SetPageRemainder
       \ifdim\PageRemainder < 97pt
        \ifdim\PageRemainder > 0pt
        \vfill\eject
       \fi\fi
    \ProcBreak \par\begingroup\parskip=0 pt%
    \goodbreak \vskip 15 pt plus 10 pt
    \noindent\null\hfill\Bibfont
      \ignorespaces #1\unskip\hfill\null\par 
    \frenchspacing \Smallfonts\rm
    \parskip=2.5 pt plus 1 pt minus.5pt%
    \nobreak\vskip 12pt plus 2pt minus2pt\nobreak
    \leftskip=0 pt \baselineskip=10.5pt}


 \def\ReferenceTagSlide{0em}
  \def\ReferenceTagGap{.5em}
  \def\ReferenceHang{30pt}

  \def \rf#1{\par\noindent
     \hangafter1\hangindent=\ReferenceHang      
     \setbox0=\hbox{\Smallfonts[\ignorespaces#1\unskip]}%
     \dimen0=\ReferenceTagGap\SetOverhang@
     \rlap{\kern\ReferenceTagSlide\box0}%
     \kern\Overhang\ignorespaces}

  \def\ref#1\par#2\par#3\par#4\par{%
     \rf{#1}#2\unskip,\ #3\unskip,\
     #4\unskip.}

  \def\endBib{\par\endgroup\vskip 12pt minus 6pt }


  \def\Coordinates{\bigskip\bigskip
     \vtop\bgroup\leftskip=\parindent
    \parskip=3pt \parindent=0pt\pretolerance=10000%
    \def\\{\hfil \break} \frenchspacing\rm }
    
  \def\endCoordinates{\par\egroup}


  \newcount\notenumber
  
  \def\note{\advance\notenumber by 1
    \footnote{\the\notenumber)}}

  \def\footnote#1{\let\@sf\empty
    \ifhmode\edef\@sf{\spacefactor\the\spacefactor}\/\fi
    \sam${}^{\fam0 #1}$\@sf\vfootnote{#1}}%

  \def\vfootnote#1{\insert\footins\bgroup
     \interlinepenalty100 \splittopskip=1pt
     \floatingpenalty=20000
     \leftskip=0pt\rightskip=0pt%
     \parindent=.3em
     \Smallfonts\rm
     \FootItem@{#1}
     \futurelet\next\fo@t}

  \def\FootItem@#1{\par\hangafter1\hangindent=\FootHang
     \setbox0=\hbox{\ignorespaces#1\unskip}%
     \dimen0=.4em\SetOverhang@
     \noindent\rlap{\box0}\kern\Overhang\ignorespaces}

  \MaxFootTag{2)}

  \def\fo@t{\ifcat\bgroup\noexpand\next \let\next\f@@t
    \else\let\next\f@t\fi \next}
  \def\f@@t{\bgroup\aftergroup\@foot\let\next}
  \def\f@t#1{\baselineskip=10pt\lineskip=1pt
            \lineskiplimit=0pt #1\@foot}%
  \def\@foot{
        \hbox{\vrule height0pt depth5pt width0pt}
        \egroup}
  \skip\footins=12 pt plus 0pt minus 0pt 
  \count\footins=1000 
  \dimen\footins=8in 


  \let\Ninepoint\ninepoint

  \def\prose{\kern1.5pt}

 \def \Blackbox
   {\leavevmode\hskip .3pt \vbox
   {\hrule height 5pt\hbox{\hskip 4.5pt}}\hskip .5pt}

 \def \XX{\Blackbox\kern.5pt\Blackbox} 

  \def\.{.\kern1pt}

    \def\Hyphen{\edef\this{\the\hyphenchar\font}%
          \hyphenchar\font=-1\char\this\hyphenchar\font=\this}

 \ifx\undefined\text
  \def\text#1{\hbox{\rm #1}}\fi 


 
 \newcount\Ht 

 \def \Acc{\expandafter }

 \def\swthat{\raise -1.1 ex\hbox{$\widehat{}$}}
 \def\swttilde{\raise -1.2 ex\hbox{$\widetilde{}$}}
 \def \overdot{{\raise .2 ex \hbox to 0pt {\hss\bf\smash{.}\hss}}}
 \def \overcircle{{\raise .1 ex \hbox to 0pt
    {\mathsurround=0pt$\eightpoint\scriptstyle\hss\circ\hss$}}}

 \def \Mathaccent#1#2{{\mathsurround=0 pt
  \setbox4=\hbox{$\vphantom{#2}$}
  \Ht=\ht4 
  \setbox5=\hbox{${#1}$}
  \setbox6=\hbox{${#2}$}
  \setbox7=\hbox to .5\wd6{}
  \copy7\kern .1\Ht \raise\Ht sp\hbox{\copy5}\kern-.1\Ht 
  \copy7\llap{\box6}
  }}

  \def\SwtCheck #1{
  \ifmmode \check{#1}%
    \else \v {#1}%
    \fi}

 \def\barpartial {%
   \kern .17 em 
    \overline {\kern -.17 em\partial\kern-.03 em}%
    \kern .03 em}

 \def\Overline#1{\setbox1=\hbox{\sam ${#1}$}%
      \ifdim \wd1 > 6pt
    \kern .11 em
    \overline {\kern -.11 em#1\kern-.14 em}
    \kern .14 em 
  \else
    \kern .03 em 
    \overline {\kern -.03 em#1\kern-.04 em}
    \kern .04 em 
  \fi}

 \def\SOverline#1{\setbox1=\hbox{\sam ${#1}$}%
      \ifdim \wd1 > 7pt
    \kern .22 em 
    \overline {\kern -.22 em#1\kern-.09 em}%
    \kern .09 em 
  \else
    \kern .10 em 
    \overline {\kern -.10 em#1\kern-.04 em}%
    \kern .04 em 
  \fi}


 \def\Underline#1{\setbox1=\hbox{\sam ${#1}$}%
      \ifdim \wd1 > 6pt
    \kern .11 em
    \underline {\kern -.11 em#1\kern-.14 em}
    \kern .14 em 
  \else
    \kern .03 em 
    \underline {\kern -.03 em#1\kern-.04 em}
    \kern .04 em 
  \fi}

 \def\SUnderline#1{\setbox1=\hbox{\sam ${#1}$}%
      \ifdim \wd1 > 7pt
    \kern .04 em 
    \underline {\kern -.04 em#1\kern-.2 em}%
    \kern .2 em 
  \else
    \kern .0 em 
    \underline {\kern -.0 em#1\kern-.15 em}%
    \kern .15 em 
  \fi}


 \ifx\MYUNDEFINED\BoxedEPSF
   \let\temp\relax
 \else
   \message{}
   \message{ !!! BoxedEPS %
         or BoxedArt macros already defined !!!}
   \let\temp\endinput
 \fi
  \temp
 
 \chardef\CatAt\the\catcode`\@
 \catcode`\@=11
 \chardef\C@tColon\the\catcode`\:
 \chardef\C@tSemicolon\the\catcode`\;
 \chardef\C@tQmark\the\catcode`\?
 \chardef\C@tEmark\the\catcode`\!

 \def\PunctOther@{\catcode`\:=12
   \catcode`\;=12 \catcode`\?=12 \catcode`\!=12}
 \PunctOther@

 \let\wlog@ld\wlog 
 \def\wlog#1{\relax} 

 \newif\ifIN@
 \newdimen\XShift@ \newdimen\YShift@ 
 \newtoks\Realtoks
 
  %
 \newdimen\Wd@ \newdimen\Ht@
 \newdimen\Wd@@ \newdimen\Ht@@
 \newdimen\TT@
 \newdimen\LT@
 \newdimen\BT@
 \newdimen\RT@
 \newdimen\XSlide@ \newdimen\YSlide@ 
 \newdimen\TheScale  
 \newdimen\FigScale  
 \newdimen\ForcedDim@@

 \newtoks\EPSFDirectorytoks@
 \newtoks\EPSFNametoks@
 \newtoks\BdBoxtoks@
 \newtoks\LLXtoks@  
 \newtoks\LLYtoks@

 \newif\ifNotIn@
 \newif\ifForcedDim@
 \newif\ifForceOn@
 \newif\ifForcedHeight@
 \newif\ifPSOrigin

 \newread\EPSFile@ 
 
  \def\ms@g{\immediate\write16}

 \newif\ifIN@\def\IN@{\expandafter\INN@\expandafter}
  \long\def\INN@0#1@#2@{\long\def\NI@##1#1##2##3\ENDNI@
    {\ifx\m@rker##2\IN@false\else\IN@true\fi}%
     \expandafter\NI@#2@@#1\m@rker\ENDNI@}
  \def\m@rker{\m@@rker}

  \newtoks\Initialtoks@  \newtoks\Terminaltoks@
  \def\SPLIT@{\expandafter\SPLITT@\expandafter}
  \def\SPLITT@0#1@#2@{\def\TTILPS@##1#1##2@{%
     \Initialtoks@{##1}\Terminaltoks@{##2}}\expandafter\TTILPS@#2@}


  \newtoks\Trimtoks@

 \def\ForeTrim@{\expandafter\ForeTrim@@\expandafter}
 \def\ForePrim@0 #1@{\Trimtoks@{#1}}
 \def\ForeTrim@@0#1@{\IN@0\m@rker. @\m@rker.#1@%
     \ifIN@\ForePrim@0#1@%
     \else\Trimtoks@\expandafter{#1}\fi}

  \def\Trim@0#1@{%
      \ForeTrim@0#1@%
      \IN@0 @\the\Trimtoks@ @%
        \ifIN@ 
             \SPLIT@0 @\the\Trimtoks@ @\Trimtoks@\Initialtoks@
             \IN@0\the\Terminaltoks@ @ @%
                 \ifIN@
                 \else \Trimtoks@ {FigNameWithSpace}%
                 \fi
        \fi
      }


   \newtoks\pt@ks
   \def \getpt@ks 0.0#1@{\pt@ks{#1}}
   \dimen0=0pt\relax\expandafter\getpt@ks\the\dimen0@

  \newtoks\Realtoks
  \def\Real#1{%
    \dimen2=#1%
      \SPLIT@0\the\pt@ks @\the\dimen2@
       \Realtoks=\Initialtoks@
            }

   \newdimen\Product
   \def\Mult#1#2{%
     \dimen4=#1\relax
     \dimen6=#2%
     \Real{\dimen4}%
     \Product=\the\Realtoks\dimen6%
        }

 \newdimen\Inverse
 \newdimen\hmxdim@ \hmxdim@=8192pt
 \def\Invert#1{%
  \Inverse=\hmxdim@
  \dimen0=#1%
  \divide\Inverse \dimen0%
  \multiply\Inverse 8}

   \def\Rescale#1#2#3{
              \divide #1 by 100\relax
              \dimen2=#3\divide\dimen2 by 100 \Invert{\dimen2}%
              \Mult{#1}{#2}%
              \Mult\Product\Inverse 
              #1=\Product}

  \def\Scale#1{\dimen0=\TheScale %
      \divide #1 by  1280 
      \divide \dimen0 by 5120 %
      \multiply#1 by \dimen0 
      \divide#1 by 10   
     }
 

 \newbox\scrunchbox

 \def\Scrunched#1{{\setbox\scrunchbox\hbox{#1}%
   \wd\scrunchbox=0pt
   \ht\scrunchbox=0pt
   \dp\scrunchbox=0pt
   \box\scrunchbox}}

 \def\Shifted@#1{%
   \vbox {\kern-\YShift@
       \hbox {\kern\XShift@\hbox{#1}\kern-\XShift@}%
           \kern\YShift@}}


 \def\cBoxedEPSF#1{{{}\leavevmode 
   \ReadNameAndScale@{#1}%
   \SetEPSFSpec@
   \ReadEPSFile@ \ReadBdB@x  
     \TrimFigDims@ 
     \CalculateFigScale@  
     \ScaleFigDims@
     \SetInkShift@
   \hbox{$\mathsurround=0pt\relax
         \vcenter{\hbox{%
             \FrameSpider{\hskip-.4pt\vrule}%
             \vbox to \Ht@{\offinterlineskip\parindent=\z@%
                \FrameSpider{\vskip-.4pt\hrule}\vfil 
                \hbox to \Wd@{\hfil}%
                \vfil
                \InkShift@{\EPSFSpecial{\EPSFSpec@}{\FigSc@leReal}}%
             \FrameSpider{\hrule\vskip-.4pt}}%
         \FrameSpider{\vrule\hskip-.4pt}}}%
     $}%
    \CleanRegisters@ 
    \ms@g{ *** Box composed for the %
         EPSF file \the\EPSFNametoks@}%
    }}      

 \def\tBoxedEPSF#1{\setbox4\hbox{\cBoxedEPSF{#1}}%
     \setbox4\hbox{\raise -\ht4 \hbox{\box4}}%
     \box4
      }

 \def\bBoxedEPSF#1{\setbox4\hbox{\cBoxedEPSF{#1}}%
     \setbox4\hbox{\raise \dp4 \hbox{\box4}}%
     \box4
      }

  \let\BoxedEPSF\cBoxedEPSF

   %

   %
  \def\gLinefigure[#1scaled#2]_#3{%
        \BoxedEPSF{#3 scaled #2}}
    
   %
   \let\EPSFfile\bBoxedEPSF
  
  \def\EPSFxsize{\afterassignment\ForceW@\ForcedDim@@}
      \def\ForceW@{\ForcedDim@true\ForcedHeight@false}
  
  \def\EPSFysize{\afterassignment\ForceH@\ForcedDim@@}
      \def\ForceH@{\ForcedDim@true\ForcedHeight@true}

  %
 \def\ReadNameAndScale@#1{\IN@0 scaled@#1@
   \ifIN@\ReadNameAndScale@@0#1@%
   \else \ReadNameAndScale@@0#1 scaled\DefaultMilScale @
   \fi}
  
 \def\ReadNameAndScale@@0#1scaled#2@{
    \let\OldBackslash@\\%
    \def\\{\OtherB@ckslash}%
    \edef\temp@{#1}%
    \Trim@0\temp@ @%
    \EPSFNametoks@\expandafter{\the\Trimtoks@ }%
    \FigScale=#2 pt%
    \let\\\OldBackslash@
    }
 
 \def\SetDefaultEPSFScale#1{%
      \global\def\DefaultMilScale{#1}}

 \SetDefaultEPSFScale{1000}

  %
 \def \SetBogusBbox@{%
     \global\BdBoxtoks@{ BoundingBox:0 0 100 100 }%
     \global\def\BdBoxLine@{ BoundingBox:0 0 100 100 }%
     \ms@g{ !!! Will use placeholder !!!}%
     }

 {\catcode`\%=12\gdef\P@S@{

 \def\ReadEPSFile@{
     \openin\EPSFile@\EPSFSpec@
     \relax  
  \ifeof\EPSFile@
     \ms@g{}%
     \ms@g{ !!! EPS FILE \the\EPSFDirectorytoks@
       \the\EPSFNametoks@\ WAS NOT FOUND !!!}
     \SetBogusBbox@
  \else
   \begingroup
   \catcode`\%=12\catcode`\:=12\catcode`\!=12
   \catcode`\G=14\catcode`\\=14\relax
   \global\read\EPSFile@ to \BdBoxLine@
   \IN@0\P@S@ @\BdBoxLine@ @%
   \ifIN@ 
     \NotIn@true
     \loop   
       \ifeof\EPSFile@\NotIn@false 
         \ms@g{}%
         \ms@g{ !!! BoundingBox NOT FOUND IN %
            \the\EPSFDirectorytoks@\the\EPSFNametoks@\ !!! }%
         \SetBogusBbox@
       \else\global\read\EPSFile@ to \BdBoxLine@
       \fi
       \global\BdBoxtoks@\expandafter{\BdBoxLine@}%
       \IN@0BoundingBox:@\the\BdBoxtoks@ @%
       \ifIN@\NotIn@false\fi%
     \ifNotIn@\repeat
   \else
         \ms@g{}%
         \ms@g{ !!! \the\EPSFNametoks@\ not PS!\  !!!}%
         \SetBogusBbox@
   \fi
  \endgroup\relax
  \fi
  \closein\EPSFile@ 
   }

  \def\ReadBdB@x{
   \expandafter\ReadBdB@x@\the\BdBoxtoks@ @}
  
  \def\ReadBdB@x@#1BoundingBox:#2@{
    \ForeTrim@0#2@%
    \IN@0atend@\the\Trimtoks@ @%
       \ifIN@\Trimtoks@={0 0 100 100 }%
         \ms@g{}%
         \ms@g{ !!! BoundingBox not found in %
         \the\EPSFDirectorytoks@\the\EPSFNametoks@\space !!!}%
         \ms@g{ !!! It must not be at end of EPSF !!!}%
         \ms@g{ !!! Will use placeholder !!!}%
       \fi
    \expandafter\ReadBdB@x@@\the\Trimtoks@ @%
   }
    
  \def\ReadBdB@x@@#1 #2 #3 #4@{
      \Wd@=#3bp\advance\Wd@ by -#1bp%
      \Ht@=#4bp\advance\Ht@ by-#2bp%
       \Wd@@=\Wd@ \Ht@@=\Ht@ 
       \LLXtoks@={#1}\LLYtoks@={#2}
      \ifPSOrigin\XShift@=-#1bp\YShift@=-#2bp\fi 
     }

   %
   \def\G@bbl@#1{}
   \bgroup
     \global\edef\OtherB@ckslash{\expandafter\G@bbl@\string\\}
   \egroup

  \def\SetEPSFDirectory{
           \bgroup\PunctOther@\relax
           \let\\\OtherB@ckslash
           \SetEPSFDirectory@}

 \def\SetEPSFDirectory@#1{
    \edef\temp@{#1}%
    \Trim@0\temp@ @
    \global\toks1\expandafter{\the\Trimtoks@ }\relax
    \egroup
    \EPSFDirectorytoks@=\toks1
    }

 \def\SetEPSFSpec@{%
     \bgroup
     \let\\=\OtherB@ckslash
     \global\edef\EPSFSpec@{%
        \the\EPSFDirectorytoks@\the\EPSFNametoks@}%
     \global\edef\EPSFSpec@{\EPSFSpec@}%
     \egroup}

  %
 \def\TrimTop#1{\advance\TT@ by #1}
 \def\TrimLeft#1{\advance\LT@ by #1}
 \def\TrimBottom#1{\advance\BT@ by #1}
 \def\TrimRight#1{\advance\RT@ by #1}

 \def\TrimBoundingBox#1{%
   \TrimTop{#1}%
   \TrimLeft{#1}%
   \TrimBottom{#1}%
   \TrimRight{#1}%
       }

 \def\TrimFigDims@{%
    \advance\Wd@ by -\LT@ 
    \advance\Wd@ by -\RT@ \RT@=\z@
    \advance\Ht@ by -\TT@ \TT@=\z@
    \advance\Ht@ by -\BT@ 
    }

  %
  \def\ForceWidth#1{\ForcedDim@true
       \ForcedDim@@#1\ForcedHeight@false}
  
  \def\ForceHeight#1{\ForcedDim@true
       \ForcedDim@@=#1\ForcedHeight@true}

  \def\ForceOn{\ForceOn@true}
  \def\ForceOff{\ForceOn@false\ForcedDim@false}
  
  \def\epsfxsize{\afterassignment\ForceW@\ForcedDim@@}
      \def\ForceW@{\ForcedDim@true\ForcedHeight@false}
  
  \def\epsfysize{\afterassignment\ForceH@\ForcedDim@@}
      \def\ForceH@{\ForcedDim@true\ForcedHeight@true}
  
  \def\CalculateFigScale@{%
     \ifForcedDim@\FigScale=1000pt
           \ifForcedHeight@
                \Rescale\FigScale\ForcedDim@@\Ht@
           \else
                \Rescale\FigScale\ForcedDim@@\Wd@
           \fi
     \fi
     \Real{\FigScale}%
     \edef\FigSc@leReal{\the\Realtoks}%
     }
   
  \def\ScaleFigDims@{\TheScale=\FigScale
      \ifForcedDim@
           \ifForcedHeight@ \Ht@=\ForcedDim@@  \Scale\Wd@
           \else \Wd@=\ForcedDim@@ \Scale\Ht@
           \fi
      \else \Scale\Wd@\Scale\Ht@        
      \fi
      \ifForceOn@\relax\else\global\ForcedDim@false\fi
      \Scale\LT@\Scale\BT@  
      \Scale\XShift@\Scale\YShift@
      }
      
 \def\HideReservedBoxes{\global\def\FrameSpider##1{\null}}
 \def\ShowReservedBoxes{\global\def\FrameSpider##1{##1}}
 \let\HideDisplacementBoxes\HideReservedBoxes  
 \let\ShowDisplacementBoxes\ShowReservedBoxes
 \let\HideFigureFrames\HideReservedBoxes
 \let\ShowFigureFrames\ShowReservedBoxes
  \ShowDisplacementBoxes
 
 \def\hSlide#1{\advance\XSlide@ by #1}
 \def\vSlide#1{\advance\YSlide@ by #1}
 
  \def\SetInkShift@{%
            \advance\XShift@ by -\LT@
            \advance\XShift@ by \XSlide@
            \advance\YShift@ by -\BT@
            \advance\YShift@ by -\YSlide@
             }
  \def\InkShift@#1{\Shifted@{\Scrunched{#1}}}
 
   %
  \def\CleanRegisters@{%
      \globaldefs=1\relax
        \XShift@=\z@\YShift@=\z@\XSlide@=\z@\YSlide@=\z@
        \TT@=\z@\LT@=\z@\BT@=\z@\RT@=\z@
      \globaldefs=0\relax}

 
 \def\SetTexturesEPSFSpecial{\PSOriginfalse
  \gdef\EPSFSpecial##1##2{\relax
    \edef\specialthis{##2}%
    \SPLIT@0.@\specialthis.@\relax
    \special{illustration ##1 scaled
                        \the\Initialtoks@}}}
 
  \def\SetUnixCoopEPSFSpecial{\PSOrigintrue 
   \gdef\EPSFSpecial##1##2{%
      \dimen4=##2pt
      \divide\dimen4 by 1000\relax
      \Real{\dimen4}
      \edef\Aux@{\the\Realtoks}%
      \includegraphics{##1\space}}}

  \def\SetBechtolsheimEPSFSpecial{\PSOrigintrue 
   \gdef\EPSFSpecial##1##2{%
      \dimen4=##2pt
      \divide\dimen4 by 1000\relax
      \Real{\dimen4}
      \edef\Aux@{\the\Realtoks}%
      \special{ps: psfiginit}%
      \special{ps: literal 1 1 0 0 1 1 startTexFig
           \the\mag\space 1000 div \Aux@\space mul 
           \the\mag\space 1000 div \Aux@\space mul scale}%
      \special{ps: include  ##1}%
      \special{ps: literal endTexFig}%
        }}

  \def\SetLisEPSFSpecial{\PSOrigintrue 
   \gdef\EPSFSpecial##1##2{%
      \dimen4=##2pt
      \divide\dimen4 by 1000\relax
      \Real{\dimen4}
      \edef\Aux@{\the\Realtoks}%
      \special{pstext="1 1 0 0 1 1 startTexFig\space
           \the\mag\space 1000 div \Aux@\space mul 
           \the\mag\space 1000 div \Aux@\space mul scale}%
      \includegraphics{##1}%
      \special{pstext=endTexFig}%
        }}

  \def\SetRokickiEPSFSpecial{\PSOrigintrue 
   \gdef\EPSFSpecial##1##2{%
      \dimen4=##2pt
      \divide\dimen4 by 10\relax
      \Real{\dimen4}
      \edef\Aux@{\the\Realtoks}%
      \includegraphics{##1}}}

  \def\SetInlineRokickiEPSFSpecial{\PSOrigintrue 
   \gdef\EPSFSpecial##1##2{%
      \dimen4=##2pt
      \divide\dimen4 by 1000\relax
      \Real{\dimen4}
      \edef\Aux@{\the\Realtoks}%
      \special{ps::[begin] 1 1 0 0 1 1 startTexFig\space
           \the\mag\space 1000 div \Aux@\space mul 
           \the\mag\space 1000 div \Aux@\space mul scale}%
      \special{ps: plotfile ##1}%
      \special{ps::[end] endTexFig}%
        }}

  \def\SetOzTeXEPSFSpecial{\PSOriginfalse 
  \gdef\EPSFSpecial##1##2{
     \special{##1\space 
       ##2 1000 div \the\mag\space 1000 div mul
       ##2 1000 div \the\mag\space 1000 div mul scale
       \the\LLXtoks@\space neg 
       \the\LLYtoks@\space neg translate
             }}} 
  
 \def\SetOzTeXPreviewedEPSFSpecial{\PSOrigintrue
 \gdef\EPSFSpecial##1##2{%
 \dimen4=##2pt
 \divide\dimen4 by 1000\relax
 \Real{\dimen4}
 \edef\Aux@{\the\Realtoks}
 \special{epsf="##1"\space scale=\Aux@}%
 }} 

  \let\SetPSprintEPSFSpecial\SetOzTeXEPSFSpecial
  \let\SetPsprintEPSFSpecial\SetOzTeXEPSFSpecial

 \def\SetArborEPSFSpecial{\PSOriginfalse 
   \gdef\EPSFSpecial##1##2{%
     \edef\specialthis{##2}%
     \SPLIT@0.@\specialthis.@\relax 
     \special{ps: epsfile ##1\space \the\Initialtoks@}}}

 \def\SetClarkEPSFSpecial{\PSOriginfalse 
   \gdef\EPSFSpecial##1##2{%
     \Rescale {\Wd@@}{##2pt}{1000pt}%
     \Rescale {\Ht@@}{##2pt}{1000pt}%
     \special{dvitops: import 
           ##1\space\the\Wd@@\space\the\Ht@@}}}

  \let\SetDVIPSONEEPSFSpecial\SetUnixCoopEPSFSpecial
  \let\SetDVIPSoneEPSFSpecial\SetUnixCoopEPSFSpecial

  \def\SetBeebeEPSFSpecial{
   \PSOriginfalse%
   \gdef\EPSFSpecial##1##2{\relax
    \special{language "PS"
      literal "##2 1000 div ##2 1000 div scale
      position = "bottom left",
      include "##1"}}}
  \let\SetDVIALWEPSFSpecial\SetBeebeEPSFSpecial

  \def\SetNorthlakeEPSFSpecial{\PSOrigintrue
   \gdef\EPSFSpecial##1##2{%
     \edef\specialthis{##2}%
     \SPLIT@0.@\specialthis.@\relax 
     \special{insert ##1,magnification=\the\Initialtoks@}}}

 \def\SetStandardEPSFSpecial{%
   \gdef\EPSFSpecial##1##2{%
     \ms@g{}
     \ms@g{%
       !!! Sorry! There is still no standard for \string%
       \special\ EPSF integration !!!}%
     \ms@g{%
      --- So you will have to identify your driver using a command}%
     \ms@g{%
      --- of the form \string\Set...EPSFSpecial, in order to get}%
     \ms@g{%
      --- your graphics to print.  See BoxedEPS.doc.}%
     \ms@g{}
     \KillEPSFSpecial
     }}

  \def\KillEPSFSpecial{\gdef\EPSFSpecial##1##2{}}

  \SetStandardEPSFSpecial 
 
 \let\wlog\wlog@ld 

 \catcode`\:=\C@tColon
 \catcode`\;=\C@tSemicolon
 \catcode`\?=\C@tQmark
 \catcode`\!=\C@tEmark

 \catcode`\@=\CatAt

 %
 %
 %
 %
 %
\let\epsffile\EPSFfile
 
  \def\FigureTitle#1{\medskip%
        \centerline{\Figurefont #1\unskip}%
        }  
    
  \catcode`\@=11
     
 \let\Mas\relax\let\mas\relax
 \let\Sam\relax\let\sam\relax
 \let\HideEdStuff\relax
 \let\PrintEdStuff\relax
 \let\KillEdStuff\relax
 \let\ShowEdStuff\relax
 \let\change\relax
 \let\beginchange\relax
 \let\endchange\relax
 \let\cbar\relax
 \let\ccbar\relax
 \let\Ninepoint\relax
 \let\ninepoint\relax
 \def\StdMathsurround{\dimen0}
 \let\LoadCMFonts\relax 
 \let\LoadPSFonts\relax 
 \let\LoadNinepoint\relax 
 \let\lBr[ 
 \let\rBr] 
 \let\SetdeGEPSFSpecial\relax
 \def\cal{\fam\tw@}

 \def\SetAuthorHead#1{}
 \def\SetTitleHead#1{}
 \def\ednote#1{}
 \def\Ednote#1{}

  \ifx\amstexloaded@\relax 
          \def\textonlyfont@#1#2{\def#1{#2}}
          \def\textfont@#1#2{\def#1{#2}}
  \else 
   \fi 
   
   \let\wlog\wlog@ld 
   \catcode`\@=\CatAt  


 

 \input amssym.def \input amssym.tex

   \def \C {\Bbb C} \def \R {\Bbb R} 
\def \M{\Bbb M}
\def \Z{\Bbb Z}
\def \C{\Bbb C}
\def \R{\Bbb R}
\def \Q{\Bbb Q}
\def \N{\Bbb N}
\def \l{\lambda}
\def \V{V^{\natural}}
\def \wt{{\rm wt}}
\def \tr{{\rm tr}}
\def \Res{{\rm Res}}
\def \End{{\rm End}}
\def \Aut{{\rm Aut}}
\def \mod{{\rm mod}}
\def \Hom{{\rm Hom}}
\def \<{\langle} 
\def \>{\rangle} 
\def \w{\omega}
\def \o{\omega}
\def \t{\tau }
\def \a{\alpha }
\def \b{\beta}
\def \e{\epsilon }
\def \L{\Lambda}
\def \la{\lambda }
\def \om{\omega }
\def \O{\Omega}
\def \voa{vertex operator algebra\ }
\def \voas{vertex operator algebras\ }
\def \p{\partial}
\def \D{\Delta}
\def \P{\Phi}

   \SetRokickiEPSFSpecial  
   \ShowFigureFrames 

 
 \H           The radical of a \voa
 
       C. Dong,\ \ H. Li,\ \ G. Mason,\ \ P. Montague

\HH 1. Introduction

Suppose that $V$ is a \voa with canonical $\Z$-grading
$$V=\coprod_{n\in\Z}V_n.$$ Each $v\in V$ has a vertex operator
$Y(v,z)=\sum_{n\in\Z}v_nz^{-n-1}$ attached to it, where $v_n\in\End
V.$ For the conformal vector $\o$ we write
$Y(\o,z)=\sum_{n\in\Z}L(n)z^{-n-2}.$ If $v$ is {\it homogeneous of
weight} $k,$ that is $v\in V_k,$ then one knows that $$v_n: V_m\to
V_{m+k-n-1}$$ and in particular the {\it zero mode} $o(v)=v_{\wt v-1}$
induces a linear operator on each $V_m.$ We extend the ``$o$'' notation
linearly to $V,$ so that in general $o(v)$ is the sum of the zero modes of the
homogeneous components of $v.$ Then we define the {\it radical} of $V$ to
be 
$$J(V)=\{v\in V|o(v)=0\}\eqno(1.1).$$

The problem arose in some work of the first three authors in [DLiM] 
and in work of the fourth author in [M] of
describing $J(V)$ precisely. We will essentially solve this problem in
the present paper in an important special case, namely that $V$ is a
\voa of {\it CFT type.} This means that the $\Z$-grading on $V$ has the shape
$$V=\coprod_{n=0}^{\infty}V_n\eqno(1.2)$$
and moreover that $V_0=\C{\bf 1}$ is spanned by the vacuum vector ${\bf 1}.$

$V$ is said to be a vertex operator algebra of {\it hermitian} CFT type if,
in addition, $V$ has the structure of a Hilbert space and further we have
an involution $v\mapsto\overline v$ on $V$ such that
$${(\overline v)_n}^\dagger=\left(e^{L(1)}v\right)_{-n}\,,$$
and $\overline\o=\o$.

In [M], the concept of a {\it deterministic} conformal field theory (a
vertex operator algebra of hermitian CFT type in our current
notation) was introduced as one for which $o(v)=0$ and $L(1)v=0$ imply
that $v\in V_1$.  Equivalently, the definition may be restated in
terms of the modes $p(v)\equiv v_{\wt v-2}$, {\it i.e.} if $L(1)v=0$
and $p(v)=0$ then $v\in V_0$. The motivation is that if we have a
state $v\in\coprod_{n=2}^\infty V_n$ such that $p(v)=0$ then
$o\left({1\over L(0)-1}v\right)$ is a generator of a (possibly
trivial) continuous symmetry of the conformal field theory.  We show
in this paper that all conformal field theories are deterministic
(restricting attention only to states which are finite sums of
components coming from distinct $V_n$'s).

To describe our result, let
$$J_1(V)=J(V)\cap V_1.\eqno(1.3)$$
We observe that $J(V)$ is in general {\it not} a $\Z$-graded subspace 
of $V,$ which accounts for some of the difficulty that its study offers, 
nevertheless the radical elements of weight one will turn out to play a special
${\rm r}\hat{\rm o}{\rm le}$. We will prove

\th Theorem 1

Suppose that $V$ is a \voa of CFT type. Then 
$$J(V)=J_1(V)+(L(0)+L(-1))V.$$
\endth

Related to Theorem 1 are the next two theorems.

\th Theorem 2

Suppose that $V$ is a \voa of CFT type. Let $v$ be homogeneous of weight at
least one. Then there is an integer $t$ such that $0\leq t\leq \wt v$ 
and
$$Y(v,z)=\sum_{n<0}v_nz^{-n-1}+\sum_{n\geq t}v_nz^{-n-1}\eqno(1.4)$$
where each operator $v_n$ in (1.4) is non-zero.
\endth

We define the {\it degree} of $v$ to be $t$ and write $\deg v=t$
if $v$ satisfies the hypotheses and conditions of Theorem 2. For
general $v\in V$ define $\deg v$ to be the least of the degrees of the
homogeneous components of $v,$ and $\deg v=-1$ if $v\in V_0.$ For
$d\geq 0$ set 
$$V^d=\{v\in V|\deg v\geq d\}.$$
Then $V^0=\coprod_{n\geq 1}V_n$ and 
$V^0\supset V^1\supset V^2\supset \cdots$ defines a filtration on
$V^0.$ We prove

\th Theorem 3

Suppose that $V$ is a \voa of CFT type. If $d\geq 1$ then
$$V^d=L(-1)^{d-1}J_1(V)+L(-1)^dV.$$
\endth

In favorable situations we will have $J_1(V)=0,$ in which case Theorems
1 and 3 take a simpler form. Although we will say something about this
situation below, we do not yet have a complete description of those
\voas of CFT type for which $J_1(V)$ is not zero. We note here
only that the \voa associated with the Heisenberg algebra has $J_1(V)\ne 0;$ 
indeed $J_1(V)=V_1$ in this case.

We refer the reader to [FHL], [FLM], [DGM] and [DHL] for
background definitions and elementary results about vertex operator algebras.

We gratefully acknowledge the following partial support:
C.D.: faculty research funds granted by the University of
California, Santa Cruz and DMS-9303374; G.M.: faculty research funds granted by the University of
California, Santa Cruz and NSF grant DMS-9401272.

\HH 2. Quasi-primary elements

Recall that $v\in V$ is called {\it quasi-primary} in case $L(1)v=0.$
It is convenient to introduce the term {\it semi-primary} for an
element $v\in V$ if $v$ is either quasi-primary or of weight 1. Note
that $v$ is quasi-primary or semi-primary if, and only if, each
homogeneous component of $v$ has the same property. 
Also note that, for a vertex operator algebra of hermitian CFT type, all
weight one states are quasi-primary and so the terms quasi- and semi-primary
are synonomous.

In this section we prove the following results under the assumption that
$V$ is a \voa of CFT type.

\th Theorem 2.1

If $v\in V$ is homogeneous and satisfies $v_{\wt v}=0,$ then $v\in V_0.$
\endth

\th Theorem 2.2

Suppose that $v\in V$ is semi-primary and satisfies $o(v)=0.$ Then $v\in V_1.$
\endth

\th Theorem 2.3

Suppose that $v\in V$ is quasi-primary and homogeneous of weight at least 2. If
$v_0=0$ then $v=0.$
\endth

The following is well-known.

\th Lemma 2.4 

The following are equivalent for a homogeneous element $v\in V:$

(a) $v\in V_0.$ 

(b) $[L(-1),v_n]=0$ for all $n\in\Z.$
\endth

\pf Proof

This follows from the identity 
$$[L(-1),Y(v,z)]=Y(L(-1)v,z)={d\over dz}Y(v,z)\eqno(2.1)$$
together with the fact that $L(-1)$ is injective on $V_n$ for $n\ne 0$
(cf. Corollary 2.4 of [DLM]). \qed

Turning to the proof of Theorem 2.1, the component version of (2.1) is
$$[L(-1),v_n]=(L(-1)v)_{n}=-nv_{n-1}.\eqno(2.2)$$
So as $v_{\wt v}=0$ then we see inductively that
$$v_n=0\eqno(2.3)$$
for $0\leq n\leq \wt v.$

Now in general one has [FLM] 
$$[a_m,b_n]=\sum_{t\geq 0}{m\choose t}(a_tb)_{m+n-t}.\eqno(2.4)$$
Taking $b_n=L(-1)=\o_0$ in (2.4) and taking (2.3) into account, we
see that
$$[v_m,L(-1)]=\sum_{t>\wt v}{m\choose t}(v_t\o)_{m-t}.\eqno(2.5)$$

Now $\wt v_t\o=\wt v-t+1,$ so $v_t\o=0$ if $t>\wt v+1.$ So (2.5)
reduces to 
$$[v_m,L(-1)]={m\choose \wt v+1}(v_{\wt v+1}\o)_{m-\wt v-1}.\eqno(2.6)$$
However, $\wt v_{\wt v+1}\o\in V_0,$ so $(v_{\wt v+1}\o)_{m-\wt v-1}=0$
unless $m-\wt v-1=-1,$ that is $m=\wt v.$ In this case ${m\choose \wt v+1}=0.$
So (2.6) yields 
$$[v_m,L(-1)]=0\eqno(2.7)$$
for all $m\in\Z.$ 

Theorem 2.1 is now a consequence of (2.7) and Lemma 2.4.
\bigbreak

We next present the proof of Theorem 2.2. Denote the homogeneous components
of $v$ by $v^i,$ that is
$$v=\sum_{i\geq 0}v^i,\ \ v^i\in V_i.\eqno(2.8)$$
Note that we do not necessarily have $o(v^i)=0,$ although each $v^i$ is 
certainly semi-primary. 

\th Lemma 2.5

The following hold:

(i) If $u$ is a homogeneous and quasi-primary then 
$$[L(1),u_t]=2(\wt u-t/2-1)u_{t+1}.\eqno(2.9)$$

(ii) If $u$ has weight 1 then
$$[L(1),u_0]=0.\eqno(2.10)$$
\endth

\pf Proof

Use (2.2) and (2.4) with $a_m=L(1)=\o_2$ to see that 
$$[L(1),u_t]=(2\wt u-t-2)u_{t+1}+(L(1)u)_t.$$
So if $u$ is quasi-primary then (2.9) holds. If $\wt u=1$ and $t=0$ then
$[L(1),u_0]=(L(1)u)_0.$ But $L(1)u\in V_0=\C{\bf 1},$ so
that $(L(1)u)_0=0.$ \qed

Turning to the proof of Theorem 2.2, since $o(v)=0$ we see from (2.8)-(2.10)
that 
$$0=[L(1),o(v)]=\sum_{i\geq 0}[L(1),o(v^i)]=\sum_{i\geq 2}2(\wt v^i-(\wt v^i-1)/2-1)v^i_{\wt v^i}.$$
So
$$\sum_{i\geq 2}(i-1)v^i_i=0.\eqno(2.11)$$

We now apply the operator $[L(1),[L(-1),*]]$ to
(2.11). We have
$$[L(1),[L(-1),u_t]]=[L(1),-tu_{t-1}]=-t(2\wt u-t-1)u_t$$
by (2.9). So (2.11)  implies that
$$\sum_{i\geq 2}(i-1)i(i-1)v^i_i=0.$$
Continuing in this fashion, we get
$$\sum_{i\geq 2}(i-1)^ki^{k-1}v^i_i=0$$
for all $k\geq 1.$ It follows that each $v^i_i=0$ for all
$i\geq 2,$ and therefore $v^i=0$ for $i\geq 2$ by Theorem 2.1. Hence
$$v=v^0+v^1.$$
To finish the theorem we need to also show that $v^0=0.$ In fact we have 
$o(v)=v^0_{-1}+v^1_0=0,$ so that $0=v^0_{-1}{\bf 1}+v^1_0{\bf 1}.$
But $v^1_0{\bf 1}=0,$ so $v^0=v^0_{-1}{\bf 1}=0,$ and the proof
of Theorem 2.2 is complete.

To prove Theorem 2.3, let $v$ be as in the statement of the theorem. Then
$v_0=0,$ and we can apply (2.9) to see inductively that $v_n=0$ for  $0\leq n\leq 2(\wt v-1).$ As $\wt v\geq 2$ then in particular
we get $v_{\wt v}=0$ and so we can apply Theorem 2.1 to 
complete the proof.

\HH 3. The shape of a vertex operator

We prove Theorems 2 and 3 in this section. Let $v\in V$ be a non-zero
homogeneous vector of weight at least one and with vertex operator
$$Y(v,z)=\sum_{n\in\Z}v_nz^{-n-1}.$$

We have $Y(v,z){\bf 1}=e^{zL(-1)}v,$ in particular if $n<0$ then $v_n{\bf 1}=
{L(-1)^{-n-1}\over (-n-1)!}v.$ As $L(-1)$ is injective on $V_n$ for
$n\ne 0$ we deduce that $v_n{\bf 1}\ne 0,$ so that $v_n\ne 0$ for $n<0.$

Next, if $v_n=0$ with $n\geq 0$ then (2.2) yields $0=[L(-1),v_n]=-nv_{n-1}.$
This then shows that $v_m=0$ for $0\leq m\leq n.$ By Theorem 2.1 we find
that $v_n\ne 0$ for $n\geq\wt v,$ and now Theorem 2 follows immediately.

Now we prove Theorem 3. We start with

\th Lemma 3.1

If $v\in V$ is homogeneous of degree $t,$ then $L(-1)^kv$ has degree $k+t$ 
for all $k\geq 0.$
\endth

\pf Proof

We may take $k=1.$  Now $Y(v,z)$ has the shape (1.4),
whence 
$$Y(L(-1)v,z)={d\over dz}Y(v,z)=-\sum_{n<0}(n+1)v_nz^{-n-2}-\sum_{n\geq t}(n+1)v_nz^{-n-2}$$
i.e.,
$$ Y(L(-1)v,z)=-\sum_{n<0}nv_{n-1}z^{-n-1}-\sum_{n\geq t+1}nv_{n-1}z^{-n-1}.\eqno(3.1)$$
Since each $v_n$ is (3.1) is non-zero we see that $\deg L(-1)v=t+1$ as
required.
\qed

\th Corollary 3.2

If $d\geq 1$ then 
$$L(-1)^{d-1}J_1(V)+L(-1)^dV\subset V^d\subset \bigoplus_{n\geq d}V_n.$$
\endth

\pf Proof

$L(-1)^{d-1}J_1(V)+L(-1)^dV\subset V^d$ follows from Lemma 3.1. The containment
$V^d\subset \oplus_{n\geq d}V_n$ reflects the fact that $\deg v\leq \wt v.$ \qed

We also need

\th Lemma 3.3

If $n\ne 1$ there is a direct sum decomposition
$$V_n=\ker (L(1):V_n\to V_{n-1})\bigoplus {\rm im}(L(-1):V_{n-1}\to V_n).$$
\endth

\pf Proof

See, for example, Proposition 3.4 of [DLM]. \qed

To Prove Theorem 3, suppose that there is $v$ of degree $d$ with 
$v\not\in L(-1)^{d-1}J_1(V)+L(-1)^dV.$ By Lemma 3.3 we can choose $v$ 
to be homogeneous and such that it has an expression of the form
$$v=\sum_{n=0}^{d-1}L(-1)^nu^n$$
with each $u^n$ homogeneous and semi-primary. Moreover we can take
$\wt v=\wt (L(-1)^nu^n),$ i.e., $\wt v=n+\wt u^n$ for $0\leq d-1.$

Now if $\wt u^n=1$ then $\wt v=n+1,$ so as $n\leq d-1$ then $\wt v\leq
d.$ So in fact $\wt v=d$ and $n=d-1$ in this case. So if $n\leq d-2$
then $u^n$ is quasi-primary of weight at least 2, so
that $\deg (L(-1)^nu^n)=n$ if $u^n\ne 0$ by Lemma 3.1 and Theorem
2.3. Since $v$ has degree $d\geq n+1$ and
$\deg(L(-1)^{d-1}u^{d-1})\geq d-1$ we see that $L(-1)^nu^n=0$ for
$n\leq d-2.$

So now $v=L(-1)^{d-1}u^{d-1}$ with $\wt u^{d-1}=1,$ forcing $\deg u^{d-1}=1$ by
Lemma 3.1. So $u^{d-1}\in J_1(V),$ and the theorem follows. \qed

\HH 4. Theorem 1

First we prove

\th Lemma 4.1 

$J_1(V)+(L(0)+L(-1))V\subset J(V).$
\endth

\pf Proof

It is enough to show that $(L(0)+L(-1))v$ lies in $J(V)$ for homogeneous 
$v\in V.$ But we have $L(0)v=(\wt v)v$ and $(L(-1)v)_n=-nv_{n-1},$ so taking 
$n=\wt (L(-1)v)-1=\wt v$ shows that $o(L(0)v+L(-1)v)=0$ as required. \qed

To begin the proof of Theorem 1, pick $v\in J(V).$ As before we can write
$$v=\sum_{n=0}^mL(-1)^nu^n\eqno(4.1)$$
where each $u^n$ is semi-primary and where $u^m\ne 0.$ We prove by induction
on $m$ that $v$ lies in  $J_1(V)+(L(0)+L(-1))V.$

Suppose first that $m=0.$ Then $v=u^m$ is semi-primary, whence the condition $o(v)=0$ forces $v\in V_1$ by Theorem 2.2. So in fact $v\in J_1(V)$ in this
case. 

In general, set $x=L(-1)^{m-1}u^m$ and $y=\sum_{n=0}^{m-1}L(-1)^nu^n.$ Thus 
$v=L(-1)x+y.$ Now from Lemma 4.1 we have $(L(0)+L(-1))x\in J(V),$ that is
$$0=o(v)=o(L(-1)x)+o(y)=-o(L(0)x)+o(y)=o(y-L(0)x).$$
We easily check that $L(0)x=(m-1)L(-1)^{m-1}u^m+L(-1)^{m-1}L(0)u^m$ so that
$$y-L(0)x=\sum_{n=0}^{m-2}L(-1)^nu^n+L(-1)^{m-1}((m-1)u^m+L(0)u^m+u^{m-1})$$
lies in $J(V).$ Since $L(0)u^m$ is semi-primary, we conclude by induction that $y-L(0)x$ lies in $J_1(V)+(L(0)+L(-1))V.$ But then the same is true of $v=y
-L(0)x+(L(0)+L(-1))x.$ This completes the proof of the theorem. 
\bigbreak

We make some remarks about Heisenberg vertex operator algebras. These are
constructed from a finite-dimensional abelian Lie algebra $H$ equipped with 
a non-degenerate, symmetric, bilinear form $\<,\>.$ One then forms the
$\Z$-graded affine Lie algebra
$$L=H\bigotimes \C[t,t^{-1}]\oplus\C c$$
where $[u\otimes t^m,v\otimes t^n]=\<u,v\>m\delta_{m+n,0}c$ and $[c,L]=0.$
The corresponding \voa has underlying Fock space
$$M(1)=S(\oplus_{n<0}H\otimes t^n)$$
with the usual action of $L.$ In particular, $c$ acts as 1.
If $h\in H$ is identified with $h\otimes
t^{-1}$ we have
$$Y(h,z)=\sum_{n\in\Z}(h\otimes t^n)z^{-n-1}$$
so that $h_n=h\otimes t^n.$ In particular, we see that $o(h)=h_0$ acts
trivially on $M(1),$ that is $o(h)=0.$ So in fact $M(1)_1=H=J_1(M(1))$ in this
case.

More generally, we see that if $M(1)$ is as above and if $V$ is any \voa then 
the tensor product $M(1)\otimes V$ (cf. [FHL]) is such that 
$M(1)\otimes {\bf 1}\subset J_1(M(1)\otimes V).$

There is a converse to this observation. To explain this, suppose that $V$ is a \voa of CFT type. On $V_1$ there is a canonical symmetric bilinear form defined by
$$\<u,v\>=u_1v.$$

\th Proposition 4.2 

Suppose that $H\subset J_1(V)$ is a subspace such that the restriction of $\<,\>$ to $H$ is non-degenerate. Then $V$ is a tensor product 
$$V\simeq M(1)\bigotimes W$$
for some \voa $W.$
\endth

\pf Proof

Given the condition on $H,$ it generates a sub \voa of $V$ (with
a different Virasoro element) isomorphic to $M(1).$ Now $V$ is a $M(1)$-module,
that is an $L$-module ($L$ as above). As is well-known (e.g., Theorem
1.7.3 of [FLM]) this implies that there is an isomorphism of $L$-modules 
$V\simeq M(1)\otimes W$ where $W$ is the space of highest weight vectors for the action of $L$ on $V.$
That is,
$$W=\{v\in V|h_nv=0, h\in H,n\geq 0\}.$$
This is the so-called commutant of $M(1)$ (cf. Theorem 5.1 of [FZ]), and
is known (loc.cit.) to be also a \voa. The proposition follows. \qed

We do not know if elements of $J_1(V)$ can occur in ways other than
those described above, though for a vertex operator algebra of {\it
hermitian} CFT type, the canonical symmetric bilinear form introduced
above is non-degenerate (see, for example, [DGM]), and so Proposition
4.2 provides the complete picture in this case.

We may consider $O_\infty(V)$, an object closely related to the
radical of $V$, which we define as follows:
$$O_\infty(V)=\{v\in V:o(v)_M=0 \ \hbox{for all modules} \ M\}\,,$$
where $o(v)_M$ is the action of the zero mode of the vertex
operator corresponding to $v$ on the module $M$.
This is to be compared to Zhu's $O(V)$, which is the set of all states
in V whose zero modes annihilate the states of lowest conformal weight
in each module. (In fact, $O_\infty(V)$ may alternatively be defined
as the intersection of objects $O_n(V)$ as described in [DLiM].)
Clearly, $O_\infty(V)\subset J(V)=J_1(V)+\left(L(0)+L(-1)\right)V$.
Further, if $v\in J_1(V)\cap O_\infty(V)$, then, by Proposition 4.2,
$V$ splits up into a tensor product $M(1)\otimes W$. We have modules
for $M(1)$ on which the zero mode of $v$ is non-zero, and so tensoring
with the adjoint module for $W$ gives a module $M$ for $V$ with
$o(v)_M\neq 0$. We deduce that $O_\infty(V)=\left(L(0)+L(-1)\right)V$.

 \Bib        Bibliography
 
 
\rf {DGM} L. Dolan, P. Goddard, P. Montague, Conformal field theories,
representations and lattice constructions, {\tt hep-th/9410029},
to appear in {\it Communications in
Mathematical Physics}.

\rf {DHL} C. Dong, H. Li, Y. Huang,
Introduction to vertex operator algebras I-III, {\it Moonshine and vertex
operator algebra,} Lecture Notes {\bf 904}, 
Research Institute for Mathematics Sciences, Kyoto University, 1995, 
1-76.

\rf {DLiM} C. Dong, H. Li, G. Mason, Vertex operator
algebras and associative algebras, preprint.

 \rf {DLM} C. Dong, Z. Lin, G. Mason,  
On vertex operator algebras as $sl_2$-modules, in: {\it Groups, 
Difference Sets, and the Monster, Proc. of a Special Research Quarter at 
The Ohio State University, Spring 1993,}  ed. by K.T. Arasu, J.F. Dillon,
K. Harada, S. Sehgal and R. Solomon, Walter de Gruyter, Berlin-New York,
1996, 349-362.

\rf {FHL} I. B. Frenkel, Y. Huang and J. Lepowsky, On
axiomatic approaches to vertex operator algebras and modules,
{\it Memoirs American Math. Soc.} {\bf 104}, 1993. 

\rf {FLM} I. B. Frenkel, J. Lepowsky, A. Meurman,
Vertex Operator Algebras and the Monster, {\it Pure and Applied
Math.,} Academic Press (1988).

\rf {FZ} I. Frenkel and Y. Zhu, Vertex operator algebras associated to
representations of affine and
Virasoro algebras, {\it Duke Math. J.} {\bf 66} (1992), 123-168.

\rf {M} P. S. Montague, Continuous symmetries of lattices and their
$\Z_2$-orbifolds, {\tt hep-th/9410218},
{\it Phys. Lett. B} {\bf 343} (1995), 113-121.

\endBib

\Coordinates
 Department of Mathematics\\
 University of California\\
 Santa Cruz, CA 95064

 Email: dong@cats.ucsc.edu, hli@crab.rutgers.edu, gem@cats.ucsc.edu, pmontagu@maths.adelaide.edu.au
 \endCoordinates
 
 \end